\documentclass[reqno]{amsart}
    \usepackage{amsmath}
    \usepackage{amssymb}
    \usepackage{amsfonts}
    \usepackage{mathrsfs}
    \usepackage[all,arc]{xy}
    \usepackage{enumerate}
    \usepackage{mathrsfs}
    \usepackage{amsthm}
    \usepackage{enumitem}
    \usepackage{extarrows}
    \usepackage{verbatim}
    \usepackage{dsfont}
     \usepackage{xcolor}
   \usepackage{smartdiagram}
    \usepackage{tensor}
    \usepackage{mathtools}
    \usepackage{float}
    \usepackage[toc,page]{appendix}
\usepackage[authoryear]{natbib}

    \theoremstyle{definition}

    \theoremstyle{remark}
    \newtheorem{rem}{Remark}

    \newcommand\smvee{\raise0.9ex\hbox{$\scriptscriptstyle\vee$}}
    \makeatletter

    \newcommand*{\rom}[1]{\expandafter\@slowromancap\romannumeral #1@}
    
    \makeatletter

\begin{document} 
\title[Comparison between CSSRBF and IMLS]{Comparison between compactly-supported spherical radial basis functions and interpolating moving least squares meshless interpolants for gravity data interpolation in geodesy and geophysics}   
\author[M. Kiani]{M. Kiani$^\dag$\\ $^\dag$\tiny Department of Geodesy, School of Surveying and Geospatial Data Engineering, University of Tehran, Tehran, Iran.}

\thanks{Corresponding author, email:mostafakiani@ut.ac.ir, tel:+989100035865.}
\date{}
\maketitle
\begin{abstract} 
The present paper is focused on the comparison of the efficiency of two important meshless interpolants for gravity acceleration interpolation. Compactly-supported spherical radial basis functions and interpolating moving least squares are used to interpolate actual gravity accelerations in southern Africa. Interpolated values are compared with actual values, gathered by observation. A thorough analysis is presented for the standard deviation of the differences between interpolated and actual values. Three different class of spherical radial basis functions-Poisson, singularity, and logarithmic-and four different type of basis functions for interpolating moving least squares approach-planar, quadratic, cubic, and spherical harmonics-are used. It is shown that in this particular problem compactly-supported spherical radial basis functions are faster and capable of achieving higher accuracies, compared to interpolating moving least squares scheme.         
\end{abstract}
$\bold{Keywords.}$ compactly-supported spherical radial basis functions, interpolating moving least squares, known points, interpolation points 
\section{Introduction}
The most important goal of the diverse class of interpolants that are used in Geodesy and Geophysics is to derive the value of the function of interest at points other than the observation points, which are located inside the region of interest, unlike the extrapolation problems (see \citet{Kiani4}). The analytical form of this function is unknown and because of this, interpolation schemes are used. There are many types of interpolant that are typically used in Geodetic problems, including spherical spline and radial basis functions (see \citet{Freeden1} and \citet{Freeden2}), ellipsoidal splines \citet{Kiani3} and spheroidal smoothing functions \citet{Kiani2}, and Interpolating Moving Least Squares (IMLS) approach (see \citet{Kiani} for more details.)

One advantage of the IMLS method over global spherical radial basis functions is that it is a local interpolant, which presents more stability and less computational complexity (refer to \citet{Wendland3}). It is logical to think of a way to use spherical radial basis functions as local interpolants. This is done using compactly-supported Spherical Radial Basis Functions (CSSRBF). 

Both the IMLS and CSSRBF are meshless, local, stable, and fast interpolants. Hence, they are very powerful schemes for real data interpolation. It is interesting to see which one is more efficient in a given problem. The comparison between the efficiencies of those two methods is the primary purpose of this paper. To test the methods, observed gravity accelerations in southern Africa are interpolated and compared with direct observations. The method that has the lowest value of standard deviation (of the difference between observed and interpolated values) is the more accurate one.

The rest of the paper is organized as follows. In section 2, a mathematical overview of the IMLS and CSSRBF methods is presented. Section 3 and 4 are devoted, respectively, to the application and the analysis of the results. In the end, conclusions are stated in section 5.

\section{Mathematical overview of IMLS and CSSRBF}
In this section, an overview of the IMLS and CSSRBF methods is presented. First, the method of IMLS is expounded. Then CSSRBF scheme is explained.

\subsection{IMLS method} This method is based on the Moving Least Squares (MLS) approach. MLS is fundamentally an approximant, which is based on Least Squares (LS) approach. The LS is performed in a local vicinity, which implies a locally supported weight function is used as the weight function for LS. Based on mathematical representations, the MLS approximant, $p_f(\xi)$, for the function $f\in C(\Omega$) ($\Omega$ is an arbitrary continuous domain), is defined as follows (see \citet{Wendland3} for more details)
\begin{equation}\label{eqn1}
p_f(\xi)=\sum_{i=1}^{n}a_i(\xi)f(\eta_i),
\end{equation}
in which $n$ is the number of points used for approximation, $a_i,~i=1,...,n$ the local coefficients of the linear approximant, and $\eta_i,~i=1,...,n$ the known points, which are the points in the vicinity of interpolation point $\xi$; they are treated as known values based on which the value at point $\xi$ is derived. The set of all known points is denoted as $D$.

The locally supported weight function, $w$, is defined based on the scaling parameter, $\delta$. So, it is needed to discuss the scaling parameter first. The scaling parameter is responsible for choosing the local, known points used for approximation. The simplest way to define the local, known points is their distance-in the sense of Euclidean norm-to the point at which the value is to be found. This point is called interpolation point hereafter. Thus, the local, known points are defined as follows

\begin{equation}\label{eqn2}
L(\xi)=\{\eta_j\in D |~||\xi-\eta_j||\le \delta\}.
\end{equation}
Now, the definition of the weight function can be presented as follows
\begin{equation}\label{eqn3}
w(\xi,\eta_j)=\psi(r),
\end{equation}
where 
\begin{equation}\label{eqn4}
\begin{split}
r&=\frac{||\xi-\eta_j||}{\delta},\\
\psi(r)&=H(1-r)\Psi(r),
\end{split}
\end{equation}
 in which $H$ is the Heaviside function, and $\Psi\in C[0,1]$ is a continuous function that represents the behavior of LS, i.e. how each known point is used in the LS procedure. Some examples of $\Psi$ are Gaussian, spline, and the reciprocal of distance (see \citet{Kiani}). After calculating all $w(\xi,\eta_i),~i=1,...,n$, the weight matrix, $W$ can be constituted as the following
\begin{equation}\label{eqn5}
W=diag[w(\xi,\eta_1),...,w(\xi,\eta_n)]. 
\end{equation} 
 
In order to calculate the local coefficients of the linear MLS approximant ($a_i,~i=1,..,n$), a set of basis functions must be chosen. In this paper, we have used four types of basis functions. These are based on the three-dimensional coordinates of the interpolation points, $\xi=(x,y,z)$, as the following
\\
\noindent$\bullet$ Planar basis functions
\begin{equation}\label{eqn6}
B_p(x,y,z)=[1,x,y,z].
\end{equation}

\noindent$\bullet$ Quadratic basis functions
\begin{equation}\label{eqn7}
B_q(x,y,z)=[1,x,y,z,x^2,y^2,z^2,xy,xz,yz].
\end{equation}

\noindent$\bullet$ Cubic basis functions
\begin{equation}\label{eqn8}
\begin{split}
b_c(x,y,z)&=[1,x,y,z,x^2,y^2,z^2,xy,xz,yz,x^3,\\
&y^3,z^3,x^2y,xy^2,x^2z,xz^2,y^2z,yz^2,xyz].
\end{split}
\end{equation}
\begin{rem}\label{rem1}
Since the three-dimensional coordinates of the interpolation points in the application presented in the next section are very large, in order to have a stable system of equations for calculating $a_i$ coefficients in \eqref{eqn1}, the parameter $\ell=\max\min||x-x_i||,~x\in\Omega,~{1\le i\le n}$ ($\eta_i=(x_i,y_i,z_i)$) is calculated and the coordinates are scaled as the following
\begin{equation}\label{eqn9}
\begin{split}
X&=\frac{x-x_i}{\ell},\\
Y&=\frac{y-y_i}{\ell},\\
Z&=\frac{z-z_i}{\ell}.
\end{split}
\end{equation}
Subsequently, the $(X,Y,Z)$ are used in the relations \eqref{eqn6}, \eqref{eqn7}, and \eqref{eqn8}.
\end{rem}
Now, the coefficients $a_i$ in \eqref{eqn1} can be computed as the following
  \begin{equation}\label{eqn10}
  \alpha=(B^TWB)^{-1}B^TWF,
  \end{equation}
in which $\alpha$ is the matrix of all the coefficients $a_i$, and $F$ is the matrix of all known values that are in $L(\xi)$.

The method of IMLS is based on interpolation and thus, the LS procedure in \eqref{eqn10} is simplified to the following system of equations
\begin{equation}\label{eqn11}
\begin{split}
p_f(\xi)&=\sum_{i=1}^{n}a_i(\xi)f(\eta_i),\\
\begin{bmatrix}
a_1\\a_2\\\vdots\\a_n
\end{bmatrix}
&=
\begin{bmatrix}
B(X_1,Y_1,Z_1)\\B(X_2,Y_2,Z_2)\\\vdots\\B(X_n,Y_n,Z_n)
\end{bmatrix}
^{-1}\times
\begin{bmatrix}
f(\eta_1)\\f(\eta_2)\\\vdots\\f(\eta_n)
\end{bmatrix}
.
\end{split}
\end{equation}
\begin{rem}\label{rem2}
The number of known points used for interpolation ($n$) is determined by the degree of basis functions. This means that based on \eqref{eqn6}-\eqref{eqn8}, $n$ is, respectively, 4, 10, 20. In general, with regard to \citet{Kiani}, the number of known points $n$ can be calculated using the following relations, based on the form of basis functions

\begin{equation}\label{eqn12}
\begin{split}
&B(x,y,z)=x^iy^jz^k,\\
&i+j+k\le v,\quad v=1,2,3,...
\end{split}
\end{equation}
\begin{equation}\label{eqn13}
\begin{split}
n&=\sum_{m=0}^{v}\binom{m+2}{2}\\
&=\frac{(v+1)(v+2)(v+3)}{6}.
\end{split}
\end{equation}
One can simply verify that for the planar, quadratic, and cubic basis function mode ($v=1$, $v=2$, and $v=3$, respectively) $n$ is 4, 10, 20.
\end{rem}
It is important to notice that the known points are chosen by the distance criterion, meaning the first $n$ points whose distances to the interpolation point is smallest are chosen as known points for interpolation in the IMLS procedure. This means the distances between interpolation point and the known points are calculated and sorted ascendingly. The first $n$ distance represent the points that must be used for interpolation. The parameter $\delta$ in \eqref{eqn2} and \eqref{eqn4} is variable for each point and is the value at which the interpolation point has exactly $n$ points in its vicinity.

Another type of basis functions used for interpolation is the spherical harmonics, up to degree $J$. It can be simply shown (see \citet{Freeden1}) that the number of spherical harmonics up to degree $J$ is $(J+1)^2$. The general form of this type of basis functions is as the following
\begin{equation}\label{eqn19}
\begin{split}
B(\phi,\lambda)&=P_{ij}(\sin\phi)\cos j\lambda,\quad i=0,...,J,~j=0,...,i,\\
B(\phi,\lambda)&=P_{ij}(\sin\phi)\sin j\lambda,\quad i=0,...,J,~j=1,...,i,
\end{split}
\end{equation}
where $P$ denotes the Legendre functions of the first kind. 
\subsection{CSSRBF method}
 This method of interpolation has such widespread use and special characteristics that the development of the theory of compactly-supported radial basis functions has been the subject of many works, including \citet{Wendland1}, \citet{Wendland2}, \citet{Schaback}, \citet{Wong}, \citet{Porcu}, and \citet{Li}. Here, the spherical radial basis functions are used as compactly-supported interpolants. This method is based on a kernel function $K$ (usually associated with the kernel of a reproducing kernel Hilbert space), which is a function of both the known and interpolation points (see \citet{Freeden1} and \citet{Freeden2} for more details.) In mathematical representation, CSSRBF method for the function $f\in C(\Omega)$ ($\Omega$ is a spherical domain), is defined as follows
 \begin{equation}\label{eqn14}
S_f(\xi)=\sum_{i=1}^{n}c_iK(\xi,\eta_i).
 \end{equation}
The coefficients $c_i,~i=1,...,n$ are determined in the same manner as those for IMLS, given in \eqref{eqn11}
\begin{equation}\label{eqn15}
\begin{split}
\begin{bmatrix}
a_1\\a_2\\\vdots\\a_n
\end{bmatrix}
&=
\begin{bmatrix}
K(\eta_1,\eta_1)&K(\eta_1,\eta_2)~~\cdots&K(\eta_1,\eta_n)\\K(\eta_2,\eta_1)&K(\eta_2,\eta_2)~~\cdots&K(\eta_2,\eta_n)\\\vdots&\vdots&\vdots\\K(\eta_n,\eta_1)&K(\eta_n,\eta_2)~~\cdots&K(\eta_n,\eta_n)
\end{bmatrix}
^{-1}\times
\begin{bmatrix}
f(\eta_1)\\f(\eta_2)\\\vdots\\f(\eta_n)
\end{bmatrix}
.
\end{split}
\end{equation} 
The role of the kernel $K$ in CSSRBF method is the same as basis functions $B$ in IMLS. In this paper, three different class of CSSRBF kernels are used. These kernels are based on the band parameter $h$ and the radius of the sphere $\Omega$, denoted by $R$, as the following

\noindent$\bullet$ Poisson kernel
\begin{equation}\label{eqn16}
K(\xi,\eta)=\frac{1}{4\pi}\frac{||\xi||^2||\eta||^2-h^2R^4}{(||\xi||^2||\eta||^2+h^2R^4-2hR^2\xi.\eta)^{\frac{3}{2}}}.
\end{equation} 
$\bullet$ Singularity kernel
\begin{equation}\label{eqn17}
K(\xi,\eta)=\frac{1}{2\pi}\frac{1}{(||\xi||^2||\eta||^2+h^2R^4-2hR^2\xi.\eta)^{\frac{1}{2}}}.
\end{equation} 
$\bullet$ Logarithmic kernel
\begin{equation}\label{eqn18}
K(\xi,\eta)=\frac{1}{4\pi R^2}\ln(1+\frac{2hR^2}{(||\xi||^2||\eta||^2+h^2R^4-2hR^2\xi.\eta)^{\frac{1}{2}}+||\xi||~||\eta||-hR^2}).
\end{equation} 
\begin{rem}\label{rem3}
The known points used for interpolation are exactly chosen in the same way as those for IMLS method.
\end{rem}
\section{Application of CSSRBF and IMLS methods in real gravity acceleration interpolation}
In this section, an application of the CSSRBF and IMLS methods is presented in gravity data interpolation in Geodesy and Geophysics. The region for interpolation is southern Africa, for which actual gravity acceleration data ($g$) are present at website www.ngdc.noaa.gov\\/mgg/gravity/. The data contain geographic latitude ($\phi$), longitude ($\lambda$), and Geodetic height ($h_G$). In order to implement the CSSRBF and IMLS methods for this problem, one must follow the diagram below

\begin{figure}[H]
	\centering
	\includegraphics[width=0.9\linewidth]{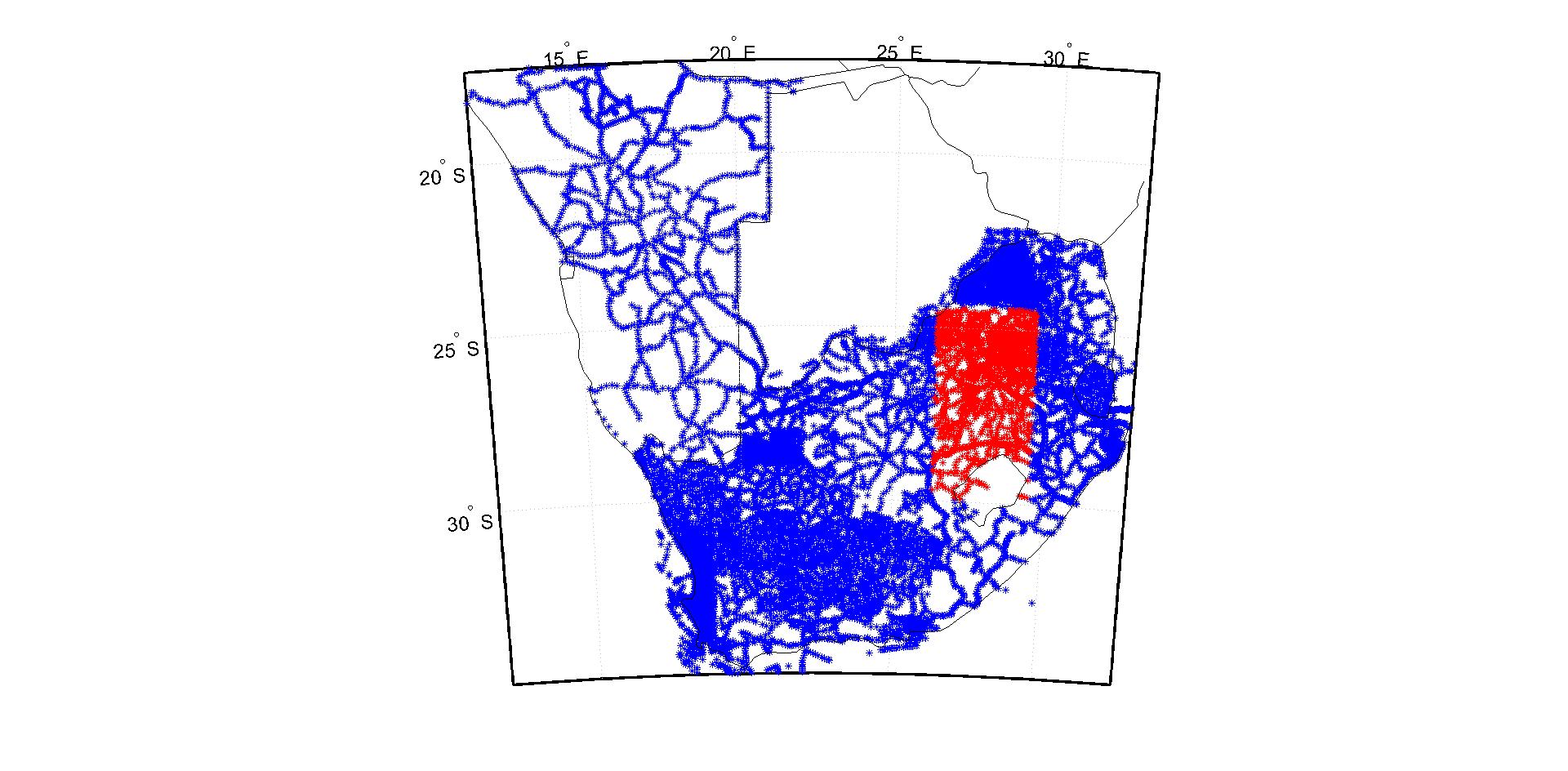}
	\caption{Known and interpolation points, respectively blue and red}
	\label{fig1}
\end{figure}

\begin{figure}[H]
	\centering
	\includegraphics[width=0.9\linewidth]{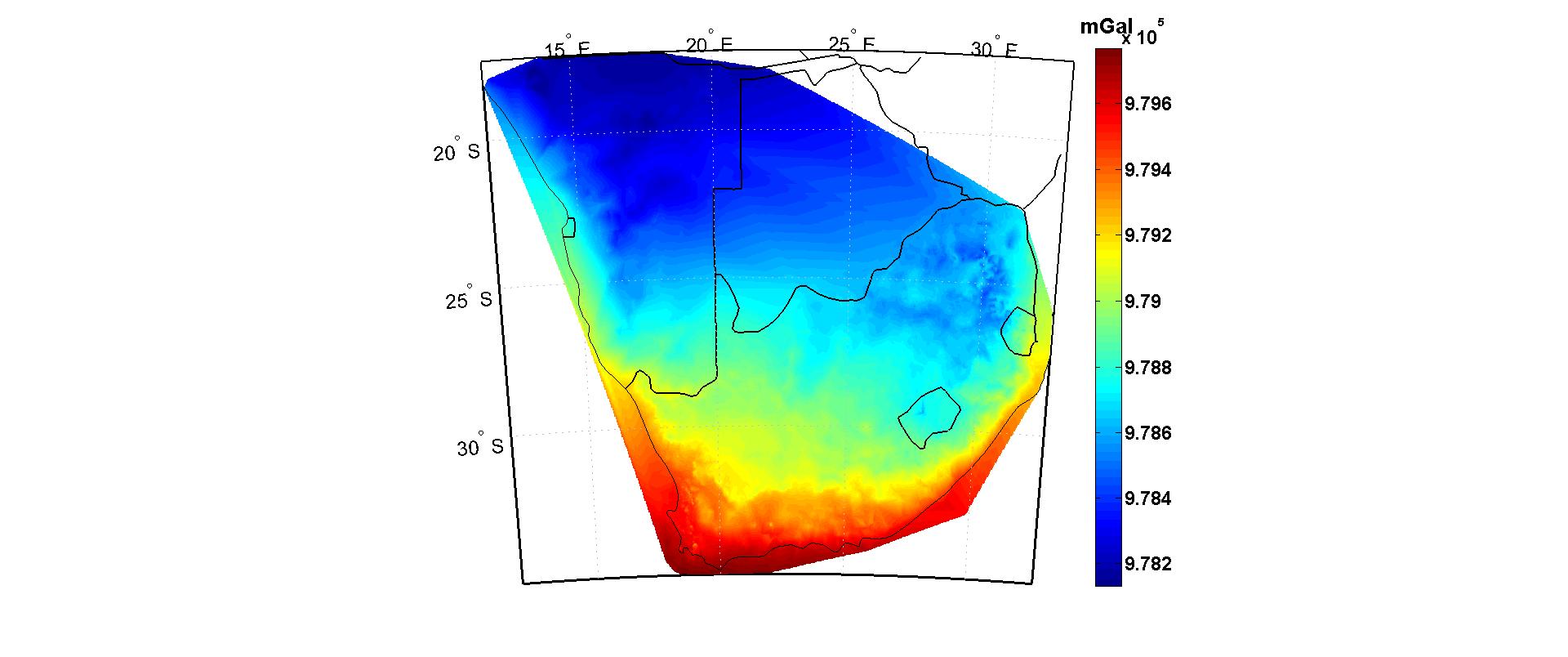}
	\caption{Observed values of all the 14559 points (known and interpolation), in mGal unit}
	\label{fig2}
\end{figure}

\begin{center}
	\smartdiagram[descriptive diagram]{
		{Step 1, {Determination of a region for interpolation points based on all the data, which contain 14559 points ($\phi,\lambda,h_G,g$)}},
		{Step 2, {Calculation of the ellipsoidal coordinates $(x,y,z)$ for both the known and interpolation points, using relations in \citet{Kiani}}},
		{Step 3, {Calculation of the distance between known and interpolation points}},
		{Step 4, {Determination of number of known points to be used for interpolation, $n$, in \eqref{eqn11} and \eqref{eqn14}}},
		{Step 5, {Caclulation of the interpolated values using the relations \eqref{eqn11} and \eqref{eqn14}}},
		{Step 6, {Subtraction of the interpolated values from observed values of interpolation points in Step 1}},
		{Step 7, {Analysis of the standard deviation ($\sigma$) of the difference of results in Step 6}}}
\end{center}
In Step 1 above, the configuration of the region of interpolation should be chosen in a way that all the interpolation points be surrounded by known points. The number of interpolation and known points is, respectively, 2000 and 12559. Figure \ref{fig1} shows the interpolation points (represented by red stars), and the known points (in blue stars).
In figure \ref{fig2}, the observed values of all 14559 points are represented.
In Table 1, the final results of Steps 3-7 are presented for planar, quadratic, and cubic basis functions.

\begin{figure}[H]
	\centering
	\includegraphics[width=0.9\linewidth]{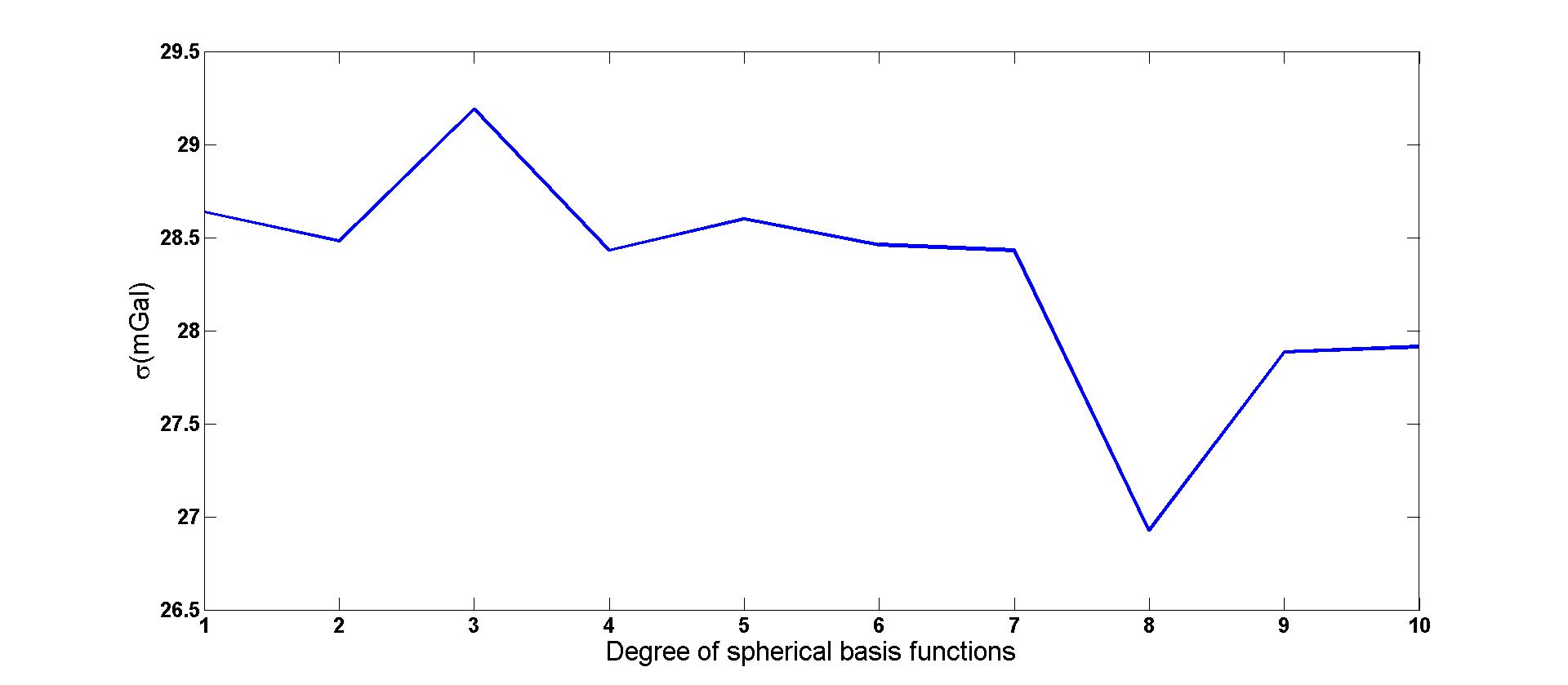}
	\caption{Standard deviation ($\sigma$) of the difference between interpolated and observed values for IMLS method with spherical harmonics basis functions}
	\label{fig12}
\end{figure}

\begin{figure}[H]
	\centering
	\includegraphics[width=0.9\linewidth]{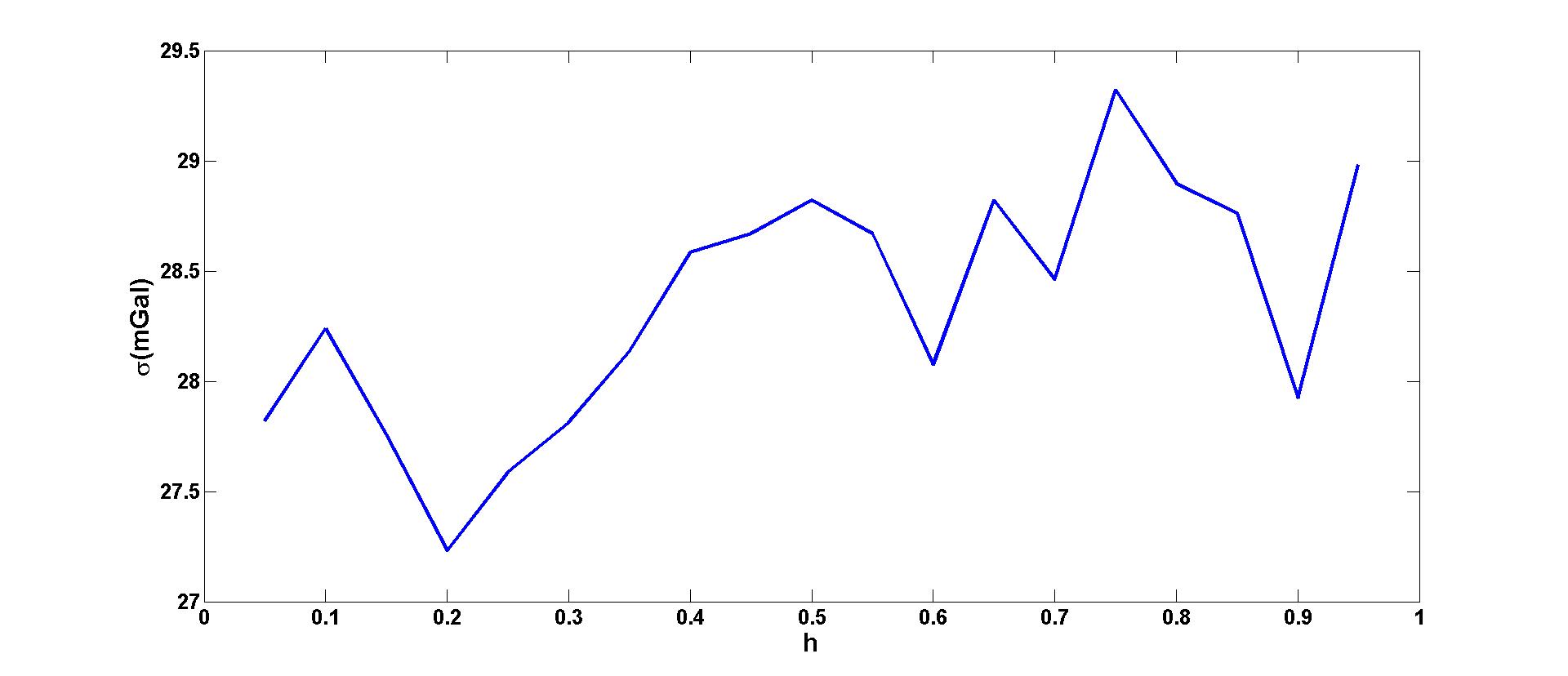}
	\caption{Standard deviation ($\sigma$) of the difference between interpolated and observed values for CSSRBF method, Poisson kernel, $n=4$}
	\label{fig3}
\end{figure}

\begin{figure}[H]
	\centering
	\includegraphics[width=0.9\linewidth]{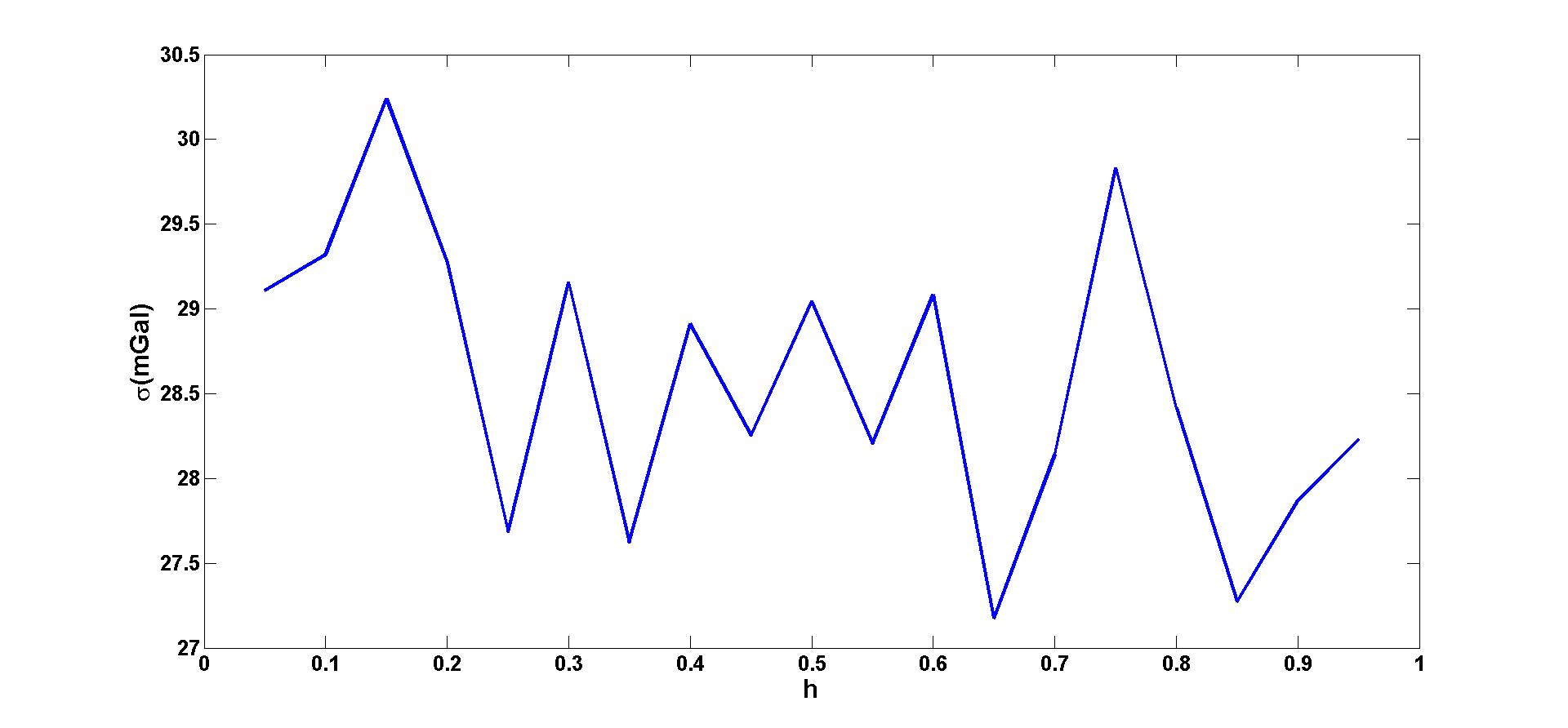}
	\caption{Standard deviation ($\sigma$) of the difference between interpolated and observed values for CSSRBF method, Poisson kernel, $n=10$}
	\label{fig4}
\end{figure}

\begin{figure}[H]
	\centering
	\includegraphics[width=0.9\linewidth]{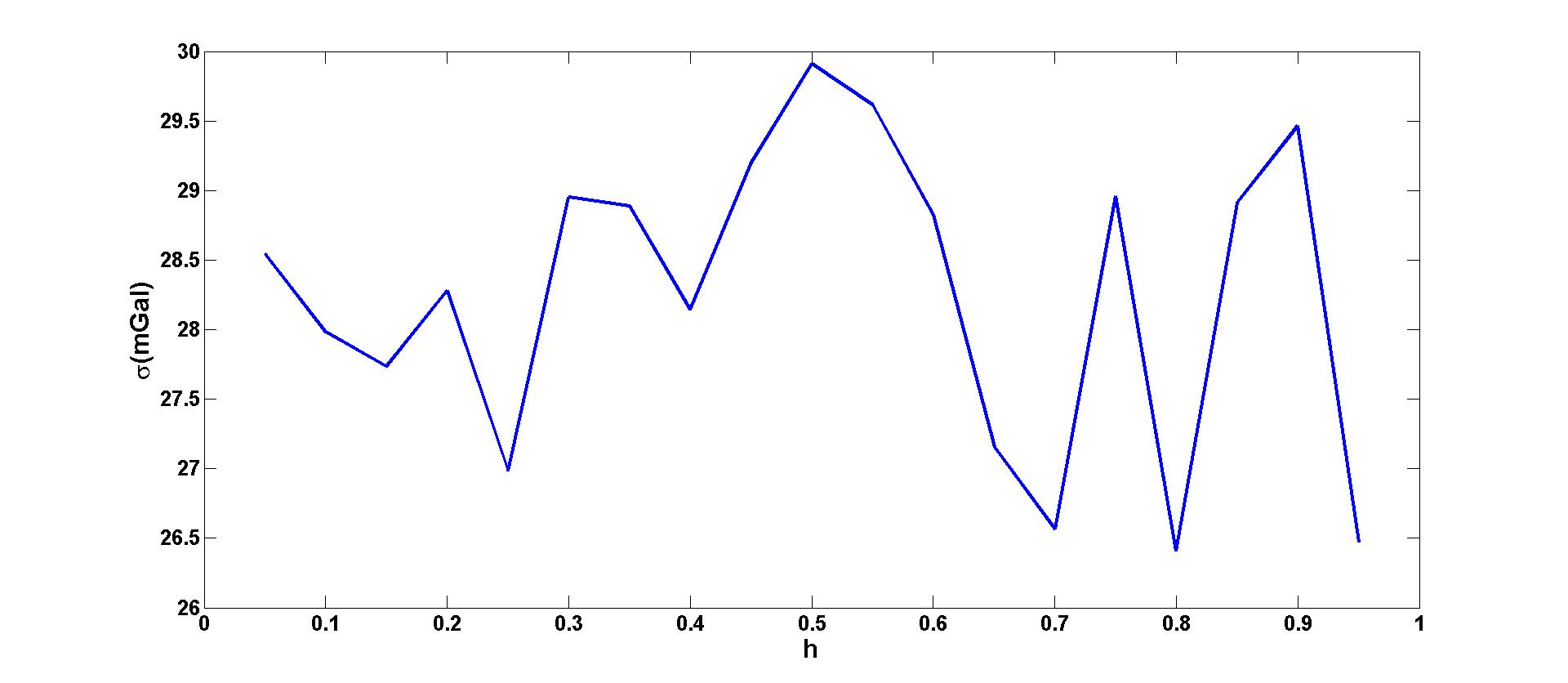}
	\caption{Standard deviation ($\sigma$) of the difference between interpolated and observed values for CSSRBF method, Poisson kernel, $n=20$}
	\label{fig5}
\end{figure}

\begin{figure}[H]
	\centering
	\includegraphics[width=0.9\linewidth]{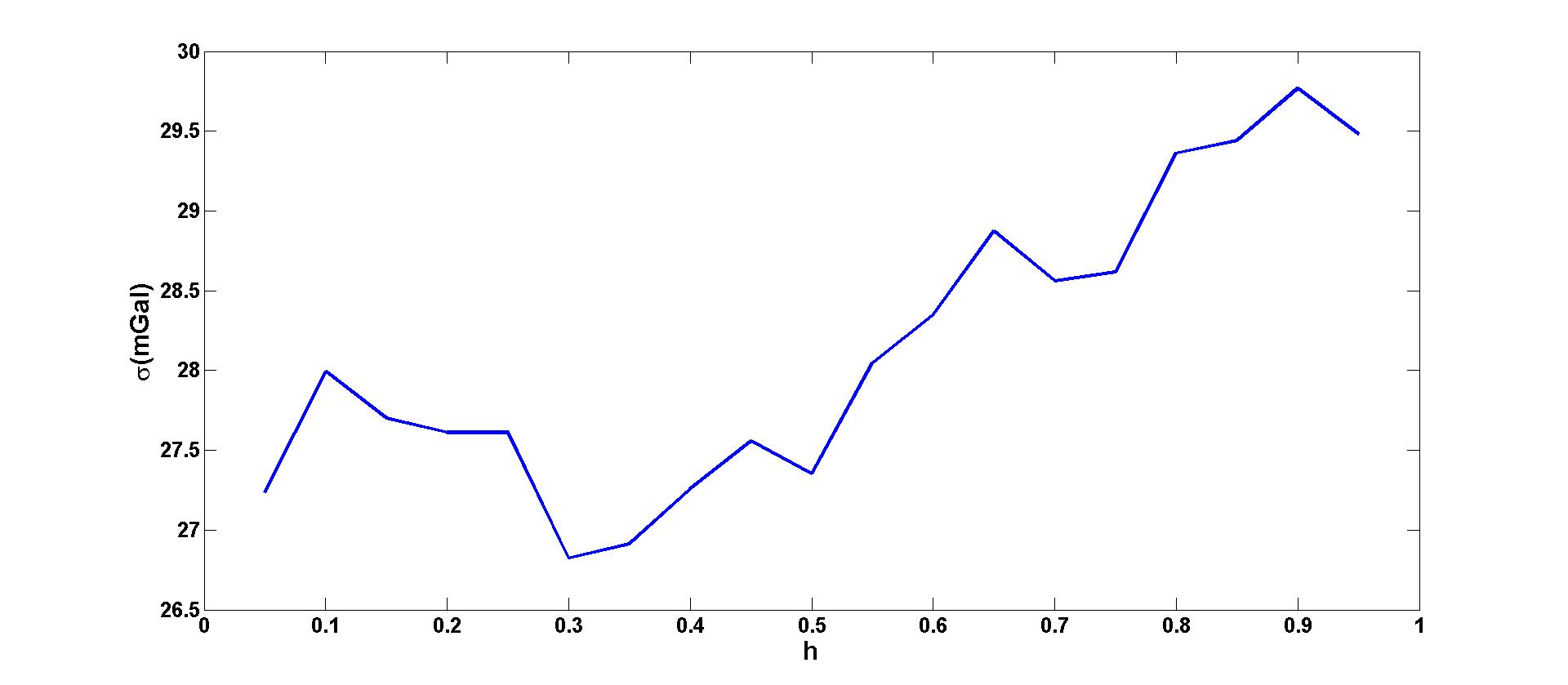}
	\caption{Standard deviation ($\sigma$) of the difference between interpolated and observed values for CSSRBF method, singularity kernel, $n=4$}
	\label{fig6}
\end{figure}

\begin{figure}[H]
	\centering
	\includegraphics[width=0.9\linewidth]{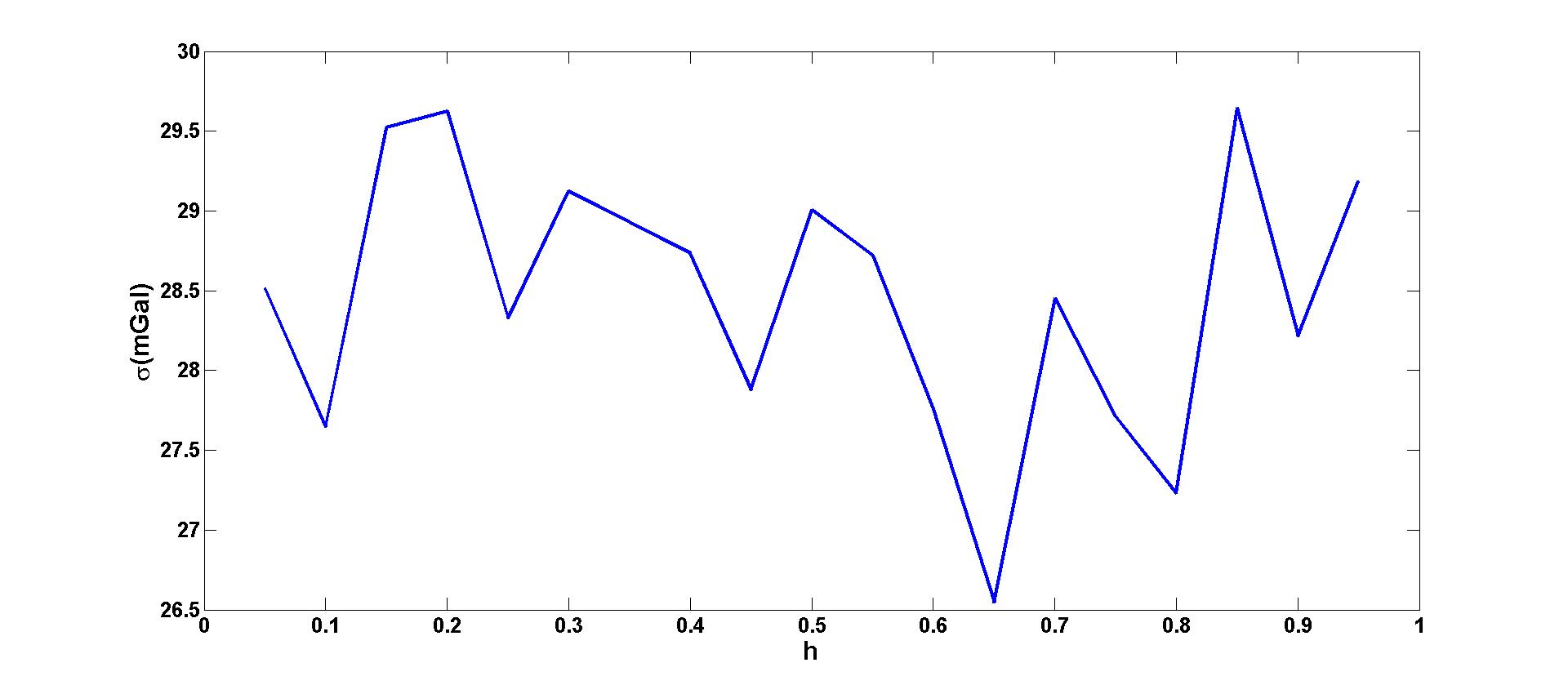}
	\caption{Standard deviation ($\sigma$) of the difference between interpolated and observed values for CSSRBF method, singularity kernel, $n=10$}
	\label{fig7}
\end{figure}

\begin{figure}[H]
	\centering
	\includegraphics[width=0.9\linewidth]{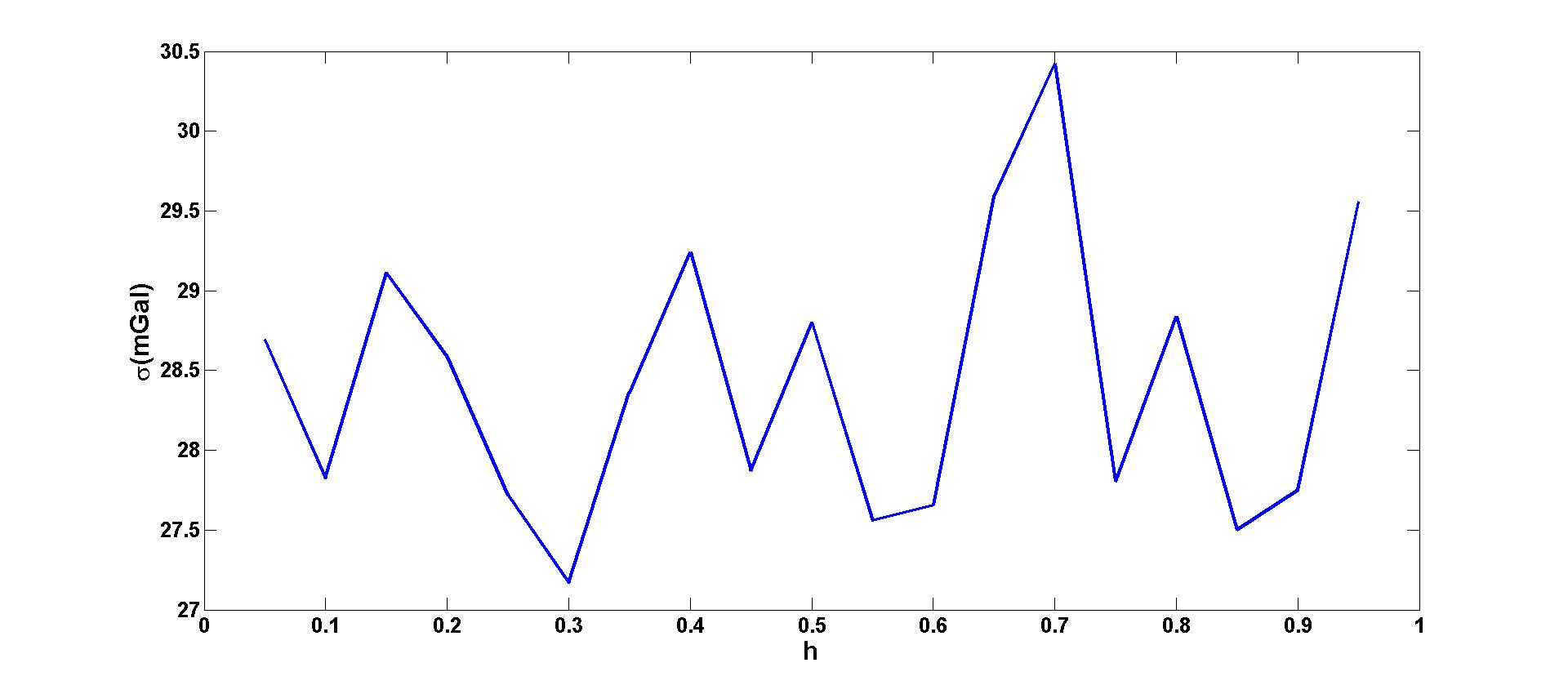}
	\caption{Standard deviation ($\sigma$) of the difference between interpolated and observed values for CSSRBF method, singularity kernel, $n=20$}
	\label{fig8}
\end{figure}

\begin{figure}[H]
	\centering
	\includegraphics[width=0.9\linewidth]{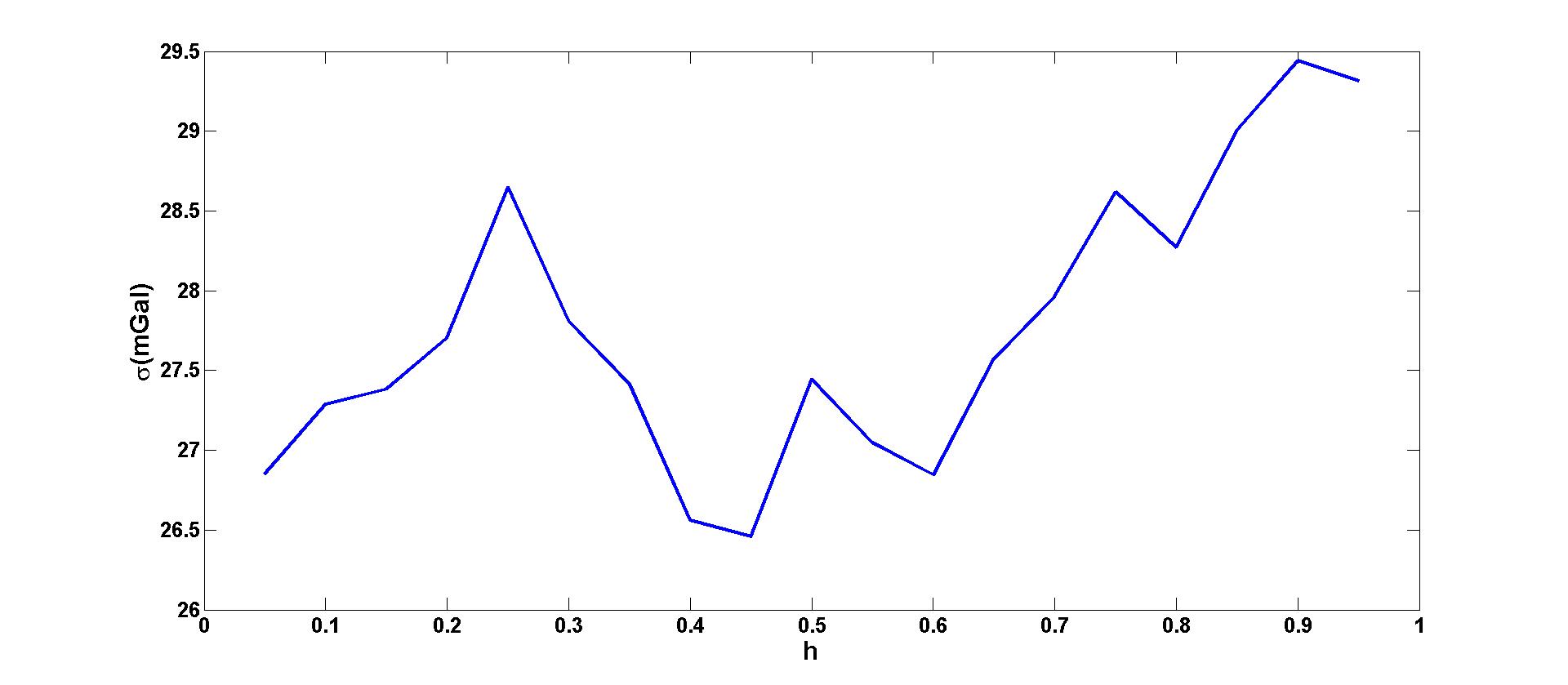}
	\caption{Standard deviation ($\sigma$) of the difference between interpolated and observed values for CSSRBF method, logarithmic kernel, $n=4$}
	\label{fig9}
\end{figure}

\begin{figure}[H]
	\centering
	\includegraphics[width=0.9\linewidth]{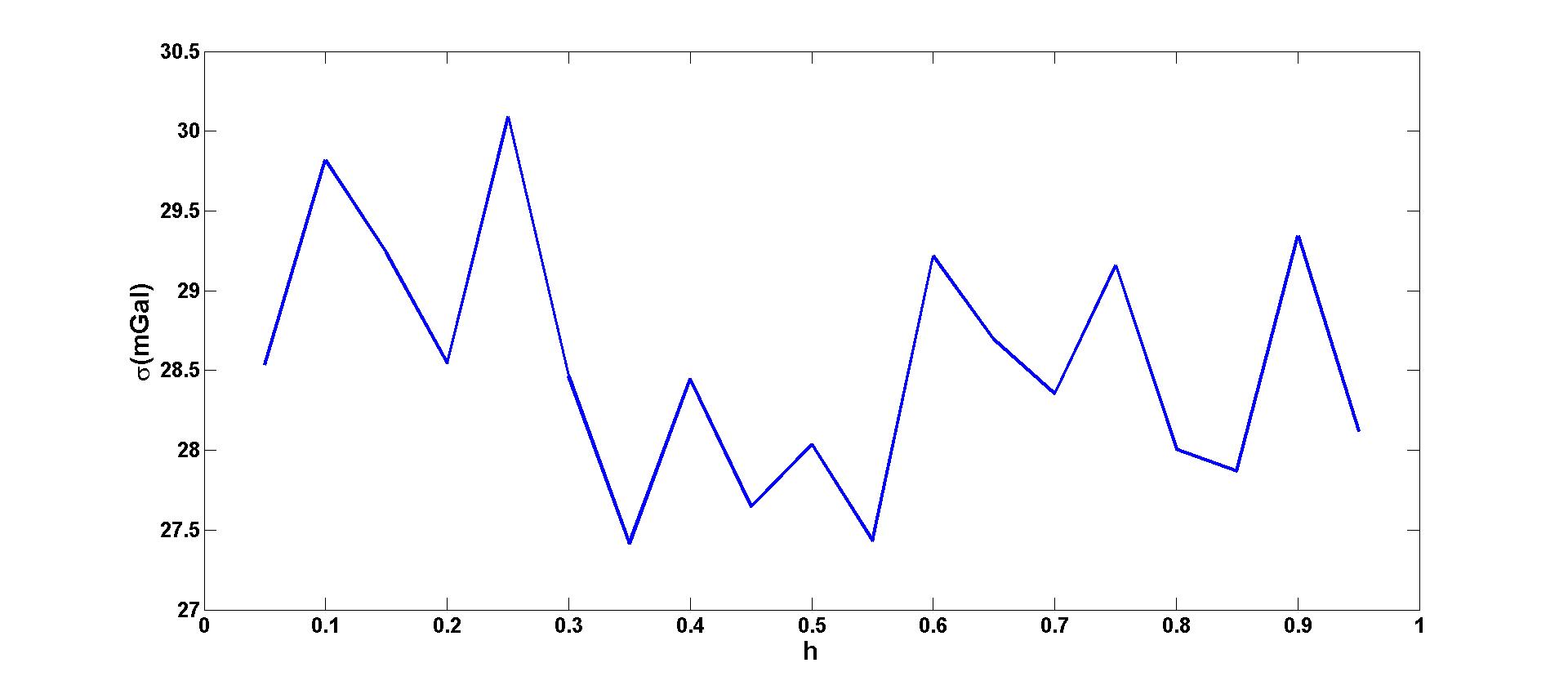}
	\caption{Standard deviation ($\sigma$) of the difference between interpolated and observed values for CSSRBF method, logarithmic kernel, $n=10$}
	\label{fig10}
\end{figure}

\begin{figure}[H]
	\centering
	\includegraphics[width=0.9\linewidth]{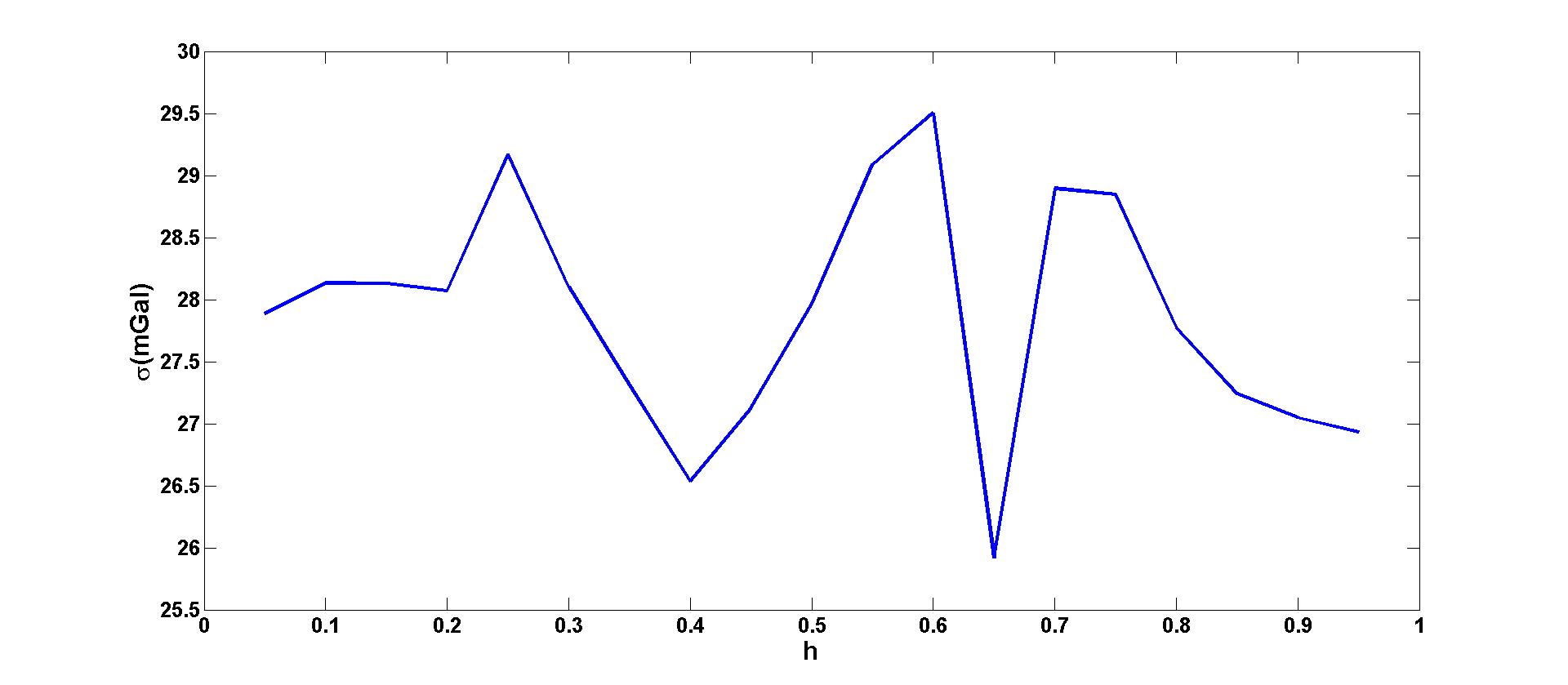}
	\caption{Standard deviation ($\sigma$) of the difference between interpolated and observed values for CSSRBF method, logarithmic kernel, $n=20$}
	\label{fig11}
\end{figure}

\begin{table}[H]
	\label{table1}
	\centering
	\caption{Standard deviation ($\sigma$) of the difference between interpolated and observed values for IMLS method}
	\begin{tabular}{|c|c|c|c|}\hline
		\cline{1-4} Type of basis functions used in IMLS procedure& Planar & Quadratic & Cubic\\
		
		\cline{1-4} $\sigma(mGal)$&27.849&28.050&27.398 \\
		\hline
	\end{tabular}
\end{table}
In figure \ref{fig12}, IMLS method with spherical basis functions is used. The degree of spherical harmonics is variable from 1-10, corresponding to 4-121 known points used for interpolation.
In the figures \ref{fig3}-\ref{fig11}, the method of CSSRBF is employed in the Steps 3-7. The three types of spherical kernel-Poisson, singularity, and logarithmic-are used for interpolation. The band parameter-$h$-in relations \eqref{eqn16}-\eqref{eqn18} is variable in each figure, ranging from 0.05 to 0.95. The number of known points used for interpolation, $n$, is also variable and is the same as those used for planar, quadratic, and cubic IMLS method.    
\section{Discussion}
In this section, an analysis of the results obtained in the previous section is presented. The important points about these results are listed below.

$\bullet$ The effect of configuration of points (i.e. how known and interpolation points are located relative to each other) is more important than the number of points used for interpolation. This can be seen in Table 1, where for planar basis functions with 4 points in the neighborhood of each interpolation point, std ($\sigma(mGal)$) is smaller than that of quadratic basis functions mode that uses 10 points for interpolation procedure. The configuration of points used for interpolation in the cubic basis functions mode, in which 20 points are used for interpolation, is better. Hence, this mode has the lowest std. 

$\bullet$ The effect of increase in the degree of spherical harmonics as basis functions on the std of difference between observed and interpolated values is slow. This can be seen from figure \ref{fig3}. The lowest value for std occurs at $v=8$, corresponding to 81 points for interpolation procedure. In this case, value of std is a balance between configuration of points and the stability of basis functions matrix in \eqref{eqn19}. As $v$ increases, the condition number of matrix of basis functions increases as well, thus making the interpolation less reliable.

$\bullet$ The lowest std value of IMLS method with spherical harmonic basis functions is approximately 26.9$~mGal$, at $v=8$; whereas the lowest std value of CSSRBF method is 25.8$~mGal$ and occurs at $h=0.8$ in the logarithmic kernel mode. It is important to notice that this std value for CSSRBF is the result of using 20 points for interpolation procedure, instead of 81 ($v=8$). This demonstrates that the CSSRBF method is more efficient than IMLS method, since the time needed for computations as well as the std value are less. In the Poisson kernel mode, the least value of std is 24.4$~mGal$, which occurs at $n=20$ and $h=0.8$. In the case of singularity kernel, values of std at $n=4$ and $h=0.3$, and $n=10$ and $h=0.65$ are smaller than that of IMLS method with spherical harmonic basis functions. Such is the case with logarithmic kernel, at $n=4$ and $h=0.45$, and $n=20$ and $h=0.65$. 

$\bullet$ At different values for $n$, the resulting std values as a function of band parameter $h$ have different behavior. In Poisson kernel mode, at $n=4$, the std is ascending; at $n=10$ descending; at $n=20$ first ascending up to $h=0.5$, and then descending. In singularity kernel mode, unlike the Poisson kernel, at all three values for $n$, the std is ascending. However, the rate of increase in the value of std is lower at $n=10$ and $n=20$, compared to $n=4$ case. This explains the role of using more points for interpolation.

$\bullet$ CSSRBF has two parameters that can be changed, namely $n$ and $h$, whereas IMLS has just the parameter $n$. This means CSSRBF is more flexible and can be used in many more different instances.      
\section{Conclusion}
A comparison between compactly-supported spherical radial basis functions and interpolating moving least squares approaches is presented in earth's gravity field data interpolation. Using both methods, at 2000 points in southern Africa the interpolated values are derived based on 12559 known points in this region. The results of interpolation are compared with actual, observed values, for different degrees of polynomial and spherical harmonic basis functions for interpolation moving least squares method, and different types of spherical kernels-Poisson, singularity, and logarithmic-for compactly-supported spherical radial basis functions scheme. It is shown that the compactly-supported spherical radial basis functions are more flexible and less time-consuming. They can be more accurate than interpolation moving least squares, as well. Hence, in this particular problem, the compactly-supported spherical radial basis functions are more efficient. However, nothing in general could be said about each individual problem. What can be concluded is that in each interpolation problem these methods should be compared to choose the best one. The comparison between compactly-supported spherical radial basis functions and interpolating moving least squares approaches with spheroidal basis functions, and in a different region of the world, could be the subject of future research, to see how these methods behave in different regions and with various types of basis functions.           
    
\end{document}